\input amssym.def
\input amssym
\input xypic


\magnification = \magstep1
\hsize=14.4truecm
\vsize= 24.5truecm
\hoffset=0.4 truecm
\baselineskip=12pt

\font\big=cmbx10 scaled \magstep2
\font\sbig=cmbx10

\font\fascio=eusm10

\font\tito=cmcsc10
\font\titol=cmcsc8
\font\proj=msbm10 scaled \magstep2
\font\eightpt=cmr8

\def\smallskip{\vskip 0.2truecm}
\def\medskip{\vskip 0.5truecm}
\def\bigskip{\vskip 0.7truecm}
\def\acapo{\par\noindent}
\def\tit{\medskip \noindent} 

\def\proof{\par\smallskip\noindent{\it Proof.}\ }
\def\endproof{\hfill$\square$\smallskip}
\definemorphism{mapsin}\solid\tip\ahook

\def\p{\Bbb P^3}
\def\pp{\Bbb P^4}
\def\ppp{\Bbb P^5}
\def\cO{{\cal O}}

\def\psf{\hbox{\fascio F}}
\def\nb#1#2{\hbox{\fascio N}_{#1 / #2}}

\def\rnb#1#2#3{\hbox{\fascio N}_{#1 / #2}\vert_{{}_ #3}}
\def\ga{p_{\scriptscriptstyle a}}
\def\cc{c_{\scriptscriptstyle 1}}

\headline={\ifnum\pageno=1\hfil\else
{\ifodd\pageno
\rightheadline\else\leftheadline\fi}\fi}
\def\rightheadline{\titol\hfil 
Threefolds in $\Bbb P^5 $ 
\dots
\hfil{\tenrm\folio}}                    
\def\leftheadline{{\tenrm\folio}\titol\hfil 
Mezzetti - Portelli \hfil}

\nopagenumbers

\vbox to 144pt{}

{\big \hbox{Threefolds in {\proj P}${}^{\hbox{\sbig 5}}$ 
with a 3-dimensional }}
\smallskip
{\big \hbox{Family of Plane Curves}}

\bigskip

\noindent{\it Emilia Mezzetti \ } and {\ \it 
Dario Portelli}
\medskip

{\eightpt {
\noindent Dipartimento di 
Scienze Matematiche, \ Universit\`a di Trieste,

\noindent Piazzale Europa 1,
 
\noindent 34127 Trieste, ITALY

\noindent e-mail: mezzette@univ.trieste.it, 
porteda@univ.trieste.it}}

\bigskip

{\eightpt 
\baselineskip=9pt
{A classification
theorem is given of smooth threefolds of $\ppp$
covered by 
a family of dimension at least three of plane 
integral curves of degree $d\geqslant 2.$  It is shown
that for such a threefold $X$ there are two possibilities:
\item{(1)} $X$ is any threefold 
contained in a hyperquadric; 

\item{(2)} $d\leqslant 3$ and $X$ is either the Bordiga or 
the Palatini scroll.
}}


\tit
{\bf Introduction}
\tit

Let $S$ be a non-degenerate surface of $\Bbb P^n$, $n\geqslant 4$. 
A classical theorem of Corrado Segre ([cS], 1921) states 
that, if $S$ contains a family of dimension $2$
of irreducible plane curves not lines, then these curves 
are conics  and $S$ is a Veronese surface of $\ppp$ or
a projection of its.  It has been remarked recently
([MP]) that this theorem is strictly related to
two  classical theorems on surfaces of $\pp$, namely the
Severi theorem and the Franchetta theorem.
The first one ([fS], 1901) asserts that a smooth surface of
 $\pp$ is 
linearly normal, unless it is a projected Veronese surface.
The theorem of Franchetta ([aF], 1947) says 
that the projected Veronese surface 
is the unique smooth surface of $\pp$ whose general
projection into $\p$ has a reducible double curve.

In the study of $3$-dimensional manifolds
 of $\ppp$ the analogous 
problems are still open. 
It has been proved that all smooth $3$-folds of $\ppp$ are 
linearly normal. 
About $2$-normality, it has been conjectured 
by Peskine and Van de Ven  that the unique non $2$-normal
smooth $3$-fold of $\ppp$ is the Palatini scroll, of degree 7
and sectional genus $4$. Moreover, when considering 
general projections of  $3$-folds of $\ppp$ into $\pp$, 
it results that the triple locus is a curve, and in all known
examples, except the Palatini scroll, this curve is irreducible.

In this paper we are interested in stating 
the analogous of Segre' s
theorem, i.e. we study $3$-dimensional manifolds of 
$\ppp$ containing a family of dimension $3$
of integral plane curves. Our results, which are contained
in Theorems 1.4 and 1.6, are the following.

Let $X$ be a smooth threefold of $\ppp$
covered by a family of dimension at least $3$
of plane integral curves not lines.
Then there are two possibilities:

\item{(1)} $X$ is any threefold contained in a hyperquadric 
($X$ is called in this case
{\it  not of isolated type});

\item{(2)} $X$ is either the Bordiga or 
the Palatini scroll and it
contains in both cases two  
families of dimension $3$ of plane curves, one
 of conics and  one of 
cubics.

\medskip

The methods that we use are mainly the adjunction theory, properties
of normal bundles and the classification of smooth 
surfaces of $\pp$ in low degree. In particular, in several situations
we heavily use the assumption that $X$ is smooth. 

\medskip

The paper is organized as follows.  
In the first section we treat the case 
of threefolds not of isolated type
and we sketch the proof of the classification theorem for those
of isolated type. In particular we recall
the results of adjunction theory that we 
shall need in the sequel.
The second section is devoted to the study of
  the normal bundle $\nb C X.$
The third section contains the rather involved proof that
a threefold of isolated type is log-special. The key point here is to
bound the self-intersection $C^2$ of $C$ on a 
hyperplane section of $X.$
In Lemma 3.1 we show that $C^2=-1$ when $C$ is a 
conic and $-1\leqslant C^2\leqslant 0$ when $C$ is a cubic.
This is the more technical part of the paper.
Finally, in sections 4 and 5 we perform 
the analysis of the log-special
threefolds not contained in a quadric and different from the 
Bordiga or the Palatini scrolls, 
to rule out the possibility that some of them contain
a family of dimension $3$ of conics or plane cubics. 

\medskip
This work has been done in the framework 
of the activities of Europroj. Both authors have been
 supported by funds of MURST and are members of GNSAGA.

\smallskip

\tit
{\bf  1.-\ Preliminaries and threefolds not of isolated type.}
\tit


Let $X\subset \Bbb P^5$ be an integral projective variety of 
dimension $3,$ and degree $d$. We will always
assume that $X$ is non-degenerate, i.e. it is not contained in 
any hyperplane. 

We suppose that $X$ contains an algebraic family {\fascio F} of 
dimension at least $3$ of plane, integral curves. 
It is immediate to remark that, if the 
union of the curves of {\fascio F} 
does not cover $X,$ then this
union is a surface containing
a 3-dimensional family of plane 
curves, hence a union of surfaces of $\p .$ 
 We will always exclude this situation,
 so from now on {\sl we assume that the curves
 of {\fascio F} cover $X.$} Our aim is to classify
such varieties $X.$

\medskip

If the curves of  {\fascio F} are lines, then 
the answer is classical and is given by the following
theorem.

\tit {\bf Theorem 1.1.} {\it Let 
$X \subset \Bbb P^5$
be an integral non-degenerate variety of dimension
 $3$
covered by a family {\fascio F} 
of dimension 
$3$ of lines. Then $X$ is ruled by planes over a curve.
If moreover $X$ is smooth, then 
$deg \ X=3$ and $X= 
\Bbb P^1\times \Bbb P^2$.}

\proof  The first assertion follows from the following theorem of 
B.Segre, the second one   from [LT].
\endproof

\tit {\bf Theorem 1.2.}(B.Segre) {\it Let $X\subset\Bbb P^n$ be a
$s$-dimensional integral projective variety, and let
$\Sigma\subset\Bbb G(k,n)$ be a component of 
maximal dimension of the variety of 
linear spaces of dimension $k$ contained in $X$.
Then:

\item{(i)} $dim\ \Sigma\leqslant (k+1)(s-k)$ and, 
if equality holds, $X$ is linear; 

\item{(ii)} if $dim\ \Sigma < (k+1)(s-k)$,
then $dim\ \Sigma\leqslant
(k+1)(s-k)-k$ and, if equality holds, then either $X$ is
a scroll in $\Bbb P^{s-1}$'s or $k=1$ and $X$ is a quadric;

\item{(iii)} if $k>1$ and $dim\ \Sigma <
(k+1)(s-k)-k$, then $dim\ \Sigma \leqslant (k+1)(s-k)-2k+1$ 
\item{{}}and, if equality holds, $X$ is a quadric.}

\proof See [bS] and [eR].
\endproof

\medskip

From now on we assume that {\sl the curves of {\fascio F}
have degree $\geqslant 2,$} so each of them generates a plane.
Under this hypothesis, we are in position to apply 
the following theorem of
classification of varieties 
containing a high dimensional 
family of degenerate subvarieties :

\tit {\bf Theorem 1.3.} {\it Let 
$X \subset \Bbb P^n$
be an integral variety of dimension
 $s$
containing a family {\fascio F} 
of dimension 
$c=h+1$ of integral subvarieties 
of dimension $s - h$. Let $Y$
be a general variety of {\fascio F} 
and assume that $Y$ spans a
$\Bbb P ^{n - h - 1}$. Then
 one of the following happens:

\item {(i)} there exists an 
integer $r$, $1\leqslant r < n - s$, such
that $X$ is contained in a
 variety $V$ of dimension at
most $n - r$ containing 
$\infty ^{h+1}$ varieties 
of dimension
$n - h - r$, each one contained 
in a linear space of dimension
$n - h - 1$;
\item {(ii)} deg $Y$ is bounded 
by a function of $h$ and $n-s$}.

\proof See [eM], Theorem (1.3).
\endproof

\smallskip 

We point out that it follows from 
the proof of the above theorem that, if $dim \ X=3,$
then possibility $(i)$ happens  precisely when 
{\sl the planes of the curves of {\fascio F} do not fill
the ambient space $\ppp .$} 

If possibility $(ii)$ happens, then $X$ is called
{\it of isolated type.} In our case, i.e. if
$n=5$, $s=3$, $h=2$, then the upper
bound on the degree of the curves of {\fascio F}
can be computed  to be $3$ ([eM], \S 3), so 
the curves  are conics or cubics.

\medskip

Next theorem gives the classification of smooth threefolds
which are not of isolated type.

 \tit {\bf Theorem 1.4} {\it Let $X$
be a smooth threefold of $\ppp $ 
of degree $d$ covered by a 
family {\fascio F} of dimension at least
$3$ of plane integral curves not lines. Then
$X$ is not of isolated type if and only if it  
is contained in a quadric. }

\proof 
Assume that $X$ is not of isolated type;
by Theorem 1.3, $X$ is contained in a variety
 $V$ of dimension $4$ containing a $3$-dimensional 
family of planes. Then $V$ is
either a quadric or a scroll in $\Bbb P^3$'s 
over a curve, by Theorem 1.2 (in the second case $V$ 
contains a family 
of planes of dimension $4$). 
Let us assume that $V$ is a scroll
 over a curve 
$B$ and let $p\colon V
\longrightarrow B$ be the map which 
gives it the scroll'\thinspace s structure.
Since $p\vert _X$ cannot be constant, 
it is surjective and the fibers
 are surfaces of $\Bbb P^3$. Hence, if 
$H$ is a general hyperplane,
 $S\colon =X\cap H$ is a smooth 
surface of $\Bbb P^4$ fibered by 
plane curves.
By a theorem of Ranestad [kR] there 
are three possibilities: 
\item {(a)} $S$ is contained in 
a quadric;
\item {(b)} $S$ is an abelian 
surface of degree 10;
\item {(c)} $S$ is a bielliptic 
surface of degree 10.
\acapo In cases $(b)$ and $(c)$ 
the irregularity $q$ of $S$ is positive: 
this  possibility is excluded by the 
assumption that $X$ is smooth, 
in view of the theorem of Barth--Lefschetz.
In case $(a)$, if $d\ > 4$  
also $X$ is contained in a quadric
by Roth's theorem (see for 
instance [MR]);
 if $d=3$ then $X=\Bbb P^1 
\times \Bbb P^2$ 
(up to projective equivalence), so 
it is contained in a quadric; 
if $d=4$ then 
$S$ must be a Del Pezzo surface,
 complete intersection of two 
quadrics, and again
$X$ is contained in a quadric.

Conversely, if $X$ is any threefold of $\Bbb P^5$
contained in a quadric $Q$, then it intersects the 
planes of $Q$ along plane curves forming an 
algebraic family of dimension at least $3$.
\endproof

\tit {\it Remarks  1.5.}  

\item {(i)} With 
notations as in 1.4, if $X$ is 
allowed to have
isolated singularities, then
$S$ is still smooth, but it 
could have $q\ >0$, 
so cases $(b)$ and $(c)$ are not excluded.

\item {(ii)} The classification of smooth 
threefolds which are contained in a quadric is complete; 
in fact
they are all arithmetically Cohen-Macaulay and either 
complete intersection or linked to a $\Bbb P^3$
in a complete intersection (see for instance [DP]).
\medskip

It remains to classify the threefolds of isolated type.

\tit {\bf Theorem 1.6} {\it Let $X$
be a smooth, non-degenerate threefold of $\ppp$ 
of degree $d,$ covered by a 
family {\fascio F} of dimension at least
$3$ of plane integral curves, not lines. Assume that
$X$ is of isolated type, i.e.
$X$ is not contained in a quadric. Then
$X$ is the Bordiga'\thinspace s scroll or
$X$ is the Palatini'\thinspace s scroll. In either cases
$X$ contains both a $3$-dimensional family of conics
and one of elliptic cubics.} 
(For the definitions of these scrolls see
e.g. [gO].)

\par\smallskip\noindent{\it Outline of the Proof.}\ 
The main tool for proving Theorem 1.6 is adjunction
theory; while giving this sketch of the proof we will
recall some fact from this theory we shall freely use 
in the sequel (for more informations see, e.g., 
[DP] and the  references given there).

Let $H$ and $K$ denote respectively a hyperplane divisor 
and a canonical divisor on $X.$ Since $X$ is not contained 
in a quadric, $X$ is different from the Segre scroll 
$\Bbb P^1\times \Bbb P^2.$ Therefore, the linear system
$\mid K+2H\mid$ is base point free and the rational map
$\phi$ associated to it, the {\it adjunction map},
is defined on the whole $X.$ Moreover, being $X$ not 
contained in a quadric, we have necessarily
$dim \ \phi (X)\geqslant 2.$

\smallskip

If $dim \ \phi (X)=2,$ then $X$ is a scroll over a 
smooth surface. Since $X$ is not contained in a quadric,
there are the following possibilities:
$X$ is the Bordiga'\thinspace s scroll or $X$
is the Palatini'\thinspace s scroll or

\smallskip

\item{(E1)} $X$ is a scroll over a $K3$ surface, 
of degree $d=9$ and sectional genus $\pi =8.$

\smallskip

If $dim \ \phi (X)=3$ then, in general, $\phi$ is the 
contraction of the exceptional planes contained in $X$.
More precisely, $\phi$ is a closed imbedding unless

\smallskip 

\item{(E2)} $X$ is of degree $d=7$ and sectional genus
$\pi =5,$ isomorphic to the 
blowing-up of a complete intersection
of three quadrics in $\Bbb P^6,$ not contained in a quadric.

\smallskip

If $X$ is such that $\phi$ is a closed imbedding, 
then the divisor $K+H$ is considered. It turns out 
that $K+H$ is nef and big unless, always under the 
assumption that $X$ is not contained in a quadric, 

\smallskip

\item{(E3)} $X$ is a Del Pezzo fibration over 
$\Bbb P^1 $ via $\mid K+H\mid ,$ with general 
fibre a complete intersection of type $(2,2)$
in $\pp .$ Moreover, $d=8$ and $\pi =7.$

\smallskip

\item{(E4)} $d=9,$\ $\pi =9,$ and 
$X$ is a conic bundle over $\Bbb P^2;$

\smallskip

\item{(E5)}  $d=12,$\ $\pi =15,$ and 
$X$ is a conic bundle over a quartic surface of $\p .$

\smallskip

\noindent
In both cases (E4), (E5) the structure of conic bundle to $X$ 
is given by the map associated to $\mid K+H\mid$.

\smallskip 

If $K+H$ is not nef and big, then $X$ is said to be of {\it 
log-special type}. 
Otherwise it is called of {\it log-general type}; 
in this last case, a general hyperplane section of $X$
is a minimal surface of general type. Moreover, a suitable
multiple of $K+H$ defines a birational morphism.

\smallskip

The key step in the proof of Theorem 1.6 is to show
that, under the given assumptions, {\sl 
$X$ must be  log-special}. As remarked after Thm. 1.3,
in the isolated case the curves of {\fascio F} can be
only conics or cubics. Then, the possibilities
(E1),$\ldots ,$ (E5) are ruled out case by case by a 
careful
analysis of the families of conics and cubic curves
lying on these threefolds. The last statement 
of the theorem is proved in [eM], \S 3. The proof is 
accomplished.

\medskip

We conclude this section with one remark.

\smallskip {\it Remark  1.7.} 
 If $X$ is of 
isolated type, 
then {\sl the family {\fascio F}
cannot have a base point.} Otherwise, 
take such a base point $P$ and project
$X$ into $\pp$ from $P.$ The image is a $3$-fold
$X'$ containing a $3$-dimensional family of lines,
hence it is a quadric or a scroll in planes, by 1.2.
Since $X$ is contained in the cone of vertex 
$P$ over $X'$, it is not of isolated type.


\bigskip\noindent
{\bf  2.-\ Some preliminaries on the normal sheaf.}
\tit 


If $C$ is a closed subvariety of a variety $X,$ 
let $\hbox{\fascio I}\subset \cO_X$ denote the 
coherent ideal which defines $C$ on $X.$ 
Recall that the {\sl normal sheaf} of $C$ on 
$X$ is defined as the dual of the $\cO_C$-module
$\hbox{\fascio I}/\hbox{\fascio I}^2\ .$ 
We shall denote this sheaf by $\nb C X \ .$

\smallskip 

It is well known that, if $C$ is locally a complete 
intersection on $X,$ e.g. if both $C$ and $X$ are 
smooth, then $\nb C X $ is a vector bundle on $C$ 
of rank equal to the codimension of $C$ in $X.$

\noindent 
In the case $C$ is an arbitrary (even singular !) 
plane curve on a smooth $X,$ then $\nb C X $ is a 
vector bundle on $C.$ In fact, the following lemma
allows us to apply the above mentioned theorem.

\smallskip

\noindent{\bf Lemma 2.1} \ {\sl 
Assume that $C\subset X\subset\Bbb P^n$ are as above. 
Assume, moreover, that $X$ is smooth and that $C$ is 
(locally) a complete intersection in $\Bbb P^n.$ 
Then $C$ is locally a complete intersection as a 
subvariety of $X\ .$}

\proof 
Being the problem local, let us consider an arbitrary 
point $x\in C\ .$  We have a surjective map\ 
$\cO_{X,x}\to \cO_{C,x}\ .$ Since $\cO_{X,x}$ is a 
regular local ring and $\cO_{C,x}$ is a
complete intersection local ring, the result follows 
from [BH], Thm. 2.3.3.
\endproof

\medskip

Let $X\subset\ppp$ be a threefold satisfying the
assumptions of Thm. 1.6. Let $C\in\hbox{\fascio F}$
be a general element and let $h\in\Bbb Q [t]$
be the Hilbert polynomial of $C.$ We will identify 
{\fascio F} with an irreducible component of the
Hilbert scheme $H$ parametrizing the
curves on $X$ having $h$ as Hilbert polynomial.
It is well known that, for $C\in H,$ the
Zariski tangent space to $H$ at $C$ is
isomorphic to $H^0(\nb C X ).$ Therefore, from
$dim (\hbox{\fascio F})\geqslant 3$ we get
$dim (H)\geqslant 3,$ hence the useful numerical
information $h^0(\nb C X )\geqslant 3$ for
any $ C\in\hbox{\fascio F}.$

\noindent
From now on, we shall denote by $\delta$
the degree of a curve $C\in\hbox{\fascio F};$
then $\delta =2$ or $\delta =3.$

\medskip

Another feature of $\nb C X $ in our situation
is given by the following

\smallskip

\noindent{\bf Lemma 2.2} \ {\sl For a
general $C\in\hbox{\fascio F}$ the bundle
$\nb C X $ is spanned by its global sections.}

\proof
Use the argument in [tF], pp. 104-105.
\endproof

\medskip

\noindent{\bf Lemma 2.3} \ {\sl Let
$dim (\hbox{\fascio F})= 3+c,$ where $c\geqslant 0.$
A general (smooth) hyperplane section $S$ of $X$ contains a 
family of dimension $c$ of curves  
$C\in\hbox{\fascio F}.$}

\proof 
Let $\diagram \varphi :\hbox{\fascio F}
\rdotted|>\tip & \Bbb G(2,5) \enddiagram$ be the rational map
which associates to the general $C\in\hbox{\fascio F}$ the
linear span $\langle C\rangle$ of $C.$ Then, the fibre of
$\varphi$ over a general element of its image has dimension $0.$
In fact, otherwise, $\langle C\rangle\subset X$ and a
straightforward application of Theorem 1.3 yields $X=\p ,$
a contradiction. Therefore, $dim (\varphi (\hbox{\fascio F}))
=3+c.$ If $\Omega\subset\Bbb G(2,5)$ represents the 
Schubert cycle parametrizing the planes of $\ppp$ contained
in a general hyperplane, then $\Omega$ intersects 
$\varphi (\hbox{\fascio F})$ properly, and this proves the
lemma since $codim(\Omega )=3.$
\endproof

\smallskip

\noindent We shall denote by $C^2$ the 
self-intersection of $C$ on $S.$ 

\smallskip

\noindent{\bf Lemma 2.4} \ {\sl We have

\noindent\item{(i)}\ $\cc(\nb C X )=C^2+\delta\ ;$

\noindent\item{(ii)}\ $K_X\cdot C=2\thinspace
\ga (C)-2-C^2-\delta\ .$}

\proof From $C\subset S \subset X$ we deduce 
the exact sequence of $\cO _C$-bundles

\vskip -0.3 cm

$$0\to\nb C S\to\nb C X\to\rnb S X C\to 0
\leqno(1)$$

\smallskip\noindent Then (i) follows,
since\  $\nb C S\simeq\cO _C(C)$\  and\  $\rnb S X C\simeq
\cO _C(1).$ 

\noindent From the exact sequence

$$0\to\hbox{\fascio T}_X\to\hbox{\fascio T}_{\ppp}\vert_{{}_X}\to
\nb X \ppp \to 0$$

\smallskip

\noindent (where $\hbox{\fascio T}_X$\ 
and\ 
$\hbox{\fascio T}_{\ppp}$ denote the tangent 
bundles to $X$ and
$\ppp$, respectively) and from 
$\cc (\hbox{\fascio T}_X)=-K_X\ $ we get

\vskip -0.3 cm

$$\cc (\nb X \ppp )=6\thinspace H\cdot X+K_X\ .
\leqno(2)$$

\noindent From the inclusions  $C\subset X
\subset \ppp$ we deduce the exact sequence of 
$\cO _C$-bundles

$$
0\to\nb C X\to\nb C \ppp\to\rnb X \ppp C\to 0
$$

\smallskip\noindent
An easy computation shows that\  $\nb C \ppp 
\simeq\cO_C(\delta )\oplus\cO_C(1)^3,$
\ hence $\cc (\nb C \ppp )=\delta^2+3\thinspace 
\delta .$  Moreover, by (2) \ 
$\cc (\rnb X \ppp C )=6\thinspace H\cdot C+K_X
\cdot C=6\thinspace \delta +K_X\cdot C,$ 
and, by the additivity of Chern classes, it 
follows 

\vskip -0.3 cm

$$K_X\cdot C=2\thinspace\ga (C)-2-\cc(\nb C X )\ .
\leqno(3)$$

\noindent 

\smallskip
\noindent Finally, combining (i) with $(3)$ we get (ii). 
\endproof

\smallskip

\noindent
For the values of $\delta$ we are interested 
in we get explicitly $K_X\cdot C=-C^2-4$ if 
$\delta =2,$ and $K_X\cdot C=-C^2-3$ for 
$\delta =3.$

\medskip


\tit
{\bf  3.-\ The proof that $X$ is log-special}
\tit

In this section we will assume that $X$ 
satisfies the hypotheses of 
Thm.\thinspace 1.6.

\medskip

\noindent{\bf Lemma 3.1}\ {\sl  On a general 
hyperplane section $S$ of $X$ we have\ $C^2=-1$
if $\delta =2,$\ and\ $C^2=0$\ or\ $C^2=-1$ if
$\delta =3.$ Moreover, $dim (\hbox{\fascio F})=3.$}
\proof We shall give the proof in several steps.

\smallskip\noindent{\tito Step 1.}\ Let $H_S$
be a general hyperplane section of $S$ containing
$C.$ Since $S$ is non-degenerate, $H_S$ is also 
non-degenerate. Therefore $H_S\neq C$ and
$H_S=C+D,$ where $D$ is a curve on $S$ of degree
$d-\delta .$ The curve $H_S$ is connected, hence 
$C\cdot D\geqslant 1.$ Assume $C\cdot D=1.$ Then
$S$ is a Veronese surface or a scroll in $\pp$
([VdV]). In the former case, $X$ is a cone over
$S.$ 
But this would imply that $X$ is not smooth, 
a contradiction. On the other hand, the smooth 
scrolls in $\pp$ are classified; they are the 
rational scroll of degree $3,$ and the 
elliptic scroll of degree $5.$ In the
former case, $X$ would be $\Bbb P^1\times\Bbb P^2,$
which is contained in a quadric hypersurface. 
The latter case is ruled out because the
irregularity of the elliptic quintic scroll 
is strictly positive, whereas we  
remarked in the proof of Thm.\thinspace 1.4, 
that the irregularity of $S$ is always zero. 
Therefore, we have $C\cdot D\geqslant 2,$ hence

$$
\delta =C\cdot H_S=C\cdot (C+D)=C^2+C\cdot D
\geqslant C^2+2$$

\smallskip

\noindent
and we conclude $C^2\leqslant 0$ for $\delta =
2,$ and $C^2\leqslant 1$ for $\delta =3.$

\smallskip\noindent{\tito Step 2.}\ We want
to exclude, now, the possibility $C^2=1$
whenever $\delta =3.$
Since $\ga (C)=1,$ the \lq\lq adjunction
formula" on $S$ yields $-C\cdot K_S=C^2.$
By Riemann-Roch 

$$
\chi (\cO_S(C))=C^2+\chi (\cO_S)=
C^2+1+h^2(\cO_S).
$$

\smallskip

\noindent
Assume $C^2=1.$ From the above relation and
from $h^2(\cO_S)\geqslant h^2(\cO_S(C))$ we get

$$
h^0(\cO_S(C))\geqslant 2+
h^2(\cO_S)-h^2(\cO_S(C))\geqslant 2
$$

\smallskip

\noindent
As a consequence, $dim \vert C\vert
\geqslant 1.$ 
If we set $h^0(\cO_S(C))= 2+k,$ 
 $k\geqslant 0,$ then it is easily checked
that $h^1(\cO_S(C))\leqslant k.$ From the long 
exact cohomology sequence associated to

$$
0\to\cO_S(C)\to\cO_S(1)\to\cO_D(1)\to 0
$$

\smallskip

\noindent we get 

$$h^0(\cO_S(C))-
h^0(\cO_S(1))+h^0(\cO_D(1))\leqslant
h^1(\cO_S(C))\leqslant k,$$

\smallskip

\noindent hence $h^0(\cO_D(H_S))\leqslant 3.$
From this inequality it follows, in particular,
that {\sl the curves $D$ are plane curves.}
Moreover, the long exact cohomology sequence 
associated to

$$
0\to\cO_S(D)\to\cO_S(1)\to\cO_C(1)\to 0
$$

\smallskip

\noindent 
jointly with $h^0(\cO_C(1))\geqslant 3$ (by
Riemann-Roch on $C$), yields 

$$h^0(\cO_S(D))\geqslant
h^0(\cO_S(1))-h^0(\cO_C(1))\geqslant 2.$$

\smallskip

\noindent
Then $\vert D\vert$ is a linear system whose
curves cover $S.$ By varing $S$ among the 
hyperplane sections of $X$ we get a family
{\fascio G} of plane curves covering $X,$
with $dim (\hbox{\fascio G})\geqslant 3.$
Since $X$ is of isolated type, we have
$deg (D)=d-3\leqslant 3,$ hence $d\leqslant 6.$
The smooth threefolds $X\subset\ppp$ with
$d<6$ are all contained in a quadric. For
$d=6$ we have that $X$ is a Bordiga'\thinspace
s scroll and $deg (D)=3.$ But the hyperplane
sections of a Bordiga'\thinspace s scroll 
contain only finitely many cubic
curves. In any case, our assumption
$C^2=1$ yields a contradiction.

\smallskip\noindent{\tito Step 3.}\ We prove, now, that
$dim(\hbox{\fascio F})=3.$

In fact, $dim(\hbox{\fascio F})
\geqslant 5$ would yield the existence on $S$ of a family of 
plane curves of dimension at least $2,$ by Lemma 2.3.
Then $S$ would be a
Veronese surface, a contradiction, since $S$ is the
hyperplane section of a smooth threefold.

If $dim(\hbox{\fascio F})=4,$ then the curves of {\fascio F}
on $S$ form a $1$-dimensional family, and 
from $C^2\leqslant 0$ it follows that this family is
actually a fibration of $S.$ 
Then, as in the proof of Thm. 1.4, 
we conclude that $X$ is contained in a quadric, 
a contradiction.

\smallskip\noindent{\tito Step 4.}\ If $C$ is a
conic, then $\nb C X\simeq
\cO_{\Bbb P^1}(a_1)\oplus\cO_{\Bbb P^1}(a_2).$ Since $\nb C X$
is spanned by its global sections, then $a_1\geqslant 0,$
$a_2\geqslant 0.$ As a consequence we have
$H^1(\nb C X)=0,$ hence, by the theory of Hilbert schemes,
{\fascio F} is smooth at any general point.
Therefore, from the previous step and from $a_1\geqslant 0,$
$a_2\geqslant 0$ we get

$$3=h^0(\nb C X )=a_1+a_2+2$$

\smallskip

\noindent and $a_1+a_2=1.$ But $a_1+a_2=c_1(\nb C X )=C^2+2,$
by Lemma 2.4. Then $C^2=-1.$

\smallskip
Therefore the proof is complete for the case
$\delta =2$ and, from now on, we will assume
$\delta =3.$

\smallskip\noindent{\tito Step 5.}\ We will prove, 
now, $C^2\geqslant -1$ in the cases
when the general $C\in\psf$ is a nodal plane
cubic or a cuspidal plane cubic (which are, a
priori, possible).
In both cases we have a  diagram

$$
\objectmargin {0.4pc}
\diagram
\Bbb P^1 \dto_n \drto^f & \cr
C  \rmapsin_j & X \cr
\enddiagram
$$

\smallskip\noindent
where $j$ is the canonical imbedding and $n$ is
the normalization map. 

\noindent For a general $C\in\psf$
by Lemma 2.2 we have that $\nb C X$ is spanned
by its global sections. Since $n^{*}$ is right
exact, $n^{*}\nb C X \simeq\cO_{\Bbb P^1}(a_1)
\oplus\cO_{\Bbb P^1}(a_2)$ is also generated by
its global sections, hence $a_1\geqslant 0$ and
$a_2\geqslant 0.$ 

\noindent Let us define $\hbox{\fascio N}_f
\ \check {} :=Ker(f^{*}\Omega_X\to\Omega_
{\Bbb P^1}).$ If $R$ denotes the ramification 
divisor on $\Bbb P^1$ of the map $n:\Bbb P^1
\to C,$ then we have an exact sequence

$$
0\to\hbox{\fascio N}_f\ \check {}\to f^{*}\Omega_X
\to\Omega_{\Bbb P^1}(-R)\to 0$$

\smallskip\noindent
([GK], proof of Lemma, pag.\thinspace 101). It
follows that $\hbox{\fascio N}_f\ \check {}$ is
a vector bundle, and we set 
$\hbox{\fascio N}_f\ \check {}\simeq
\cO_{\Bbb P^1}(-b_1)\oplus\cO_{\Bbb P^1}(-b_2).$
Moreover, the degree of $\hbox{\fascio N}_f\ 
\check {}$ or, better, of $\hbox{\fascio N}_f$
is easily computed to be

$$
b_1+b_2=C^2+1\ \ \ \ \ \ \ \hbox{for the nodal case;}
$$

\vskip -0.9 cm

$$\leqno(4)$$

\vskip -0.9 cm

$$
b_1+b_2=C^2\ \ \ \ \ \ \  \ \ 
\hbox{for the cuspidal case.}
$$

\smallskip
\noindent
Now, we have the canonical, commutative
diagram with exact rows

$$
\objectmargin {0.4pc}
\diagram
& n^{*}\nb C X \ \check {} \rto & f^{*}\Omega_X 
\rto \ddouble & n^{*}\Omega_C \dto \rto & 0\cr
0 \rto & \hbox{\fascio N}_f\ \check {} \rto &
f^{*}\Omega_X \rto & \Omega_{\Bbb P^1} & \cr
\enddiagram\leqno(5)
$$

\smallskip\noindent
Then there exists a map $\alpha :
 n^{*}\nb C X \ \check {}\to\hbox{\fascio N}_f\ 
\check {} .$ This map is represented by a matrix

$$\left(\matrix{F_{11} & F_{12}\cr
F_{21} & F_{22} \cr}\right)$$

\smallskip\noindent where the $F_{ij}$'s are
homogeneous polynomials of degree $a_i-b_j.$ For
future use we want to determine the degree of the
determinant of $\alpha .$ Of course, it is equal
to $dim_k(Coker(\alpha ));$ moreover, it is
easily seen that the support of $(Coker(\alpha ))$ 
is contained into $n^{-1}(Sing(C)).$ Finally, let 
us consider the exact sequence

$$
0\to K\to n^{*}\Omega_C \to\Omega_{\Bbb P^1}\to
\Omega_{\Bbb P^1/C}\to 0$$

\smallskip\noindent
Let $x\in n^{-1}(Sing(C));$ by applying the \lq\lq
snake'\thinspace s lemma" to the diagram obtained 
localizing $(6)$ at $x ,$ we show that $Coker(
\alpha )_x\simeq K_x\ .$ Now, in the nodal case
$\Omega_{\Bbb P^1/C}=0.$ Therefore $dim_k(K_x)=1.$
Since $n^{-1}(Sing(C))$ consists of two points, we
get $dim_k(Coker(\alpha ))=2.$ In the cuspidal case
$\Omega_{\Bbb P^1/C}$ is a skyscraper sheaf on the
unique point $x\in n^{-1}(Sing(C)),$ and
$dim_k(\Omega_{\Bbb P^1/C,x})=1.$ Therefore,
$K$ is a skyscraper sheaf on $x,$ and
$dim_k(K)=2.$ Hence, both in the nodal and in the
cuspidal case we have

\vskip -0.3 cm

$$
a_1-b_2+a_2-b_1=deg(det(\alpha ))=2.
\leqno(6)$$


After these preparations we will deal first with the
cuspidal case. If we tensorize the canonical inclusion
$\cO_C\subset n_{*}\cO_{\Bbb P^1}$ by $\nb C X ,$
then we get by
the \lq\lq projection formula" another inclusion, namely
$\nb C X \subset n_{*}n^{*}\nb C X .$ Therefore

$$
3\leqslant h^0(\nb C X )\leqslant h^0(n^{*}\nb C X )
=a_1+a_2+2 ,
$$

\smallskip\noindent
$a_1+a_2\geqslant 1.$ From $(6)$ we get
$b_1+b_2\geqslant -1,$ and the conclusion follows
by $  (4).$

\smallskip

Rather surprisingly, this approach is too rough for
the nodal case; in fact, it yields only $C^2
\geqslant -2.$ To overcome this difficulty, we intrduce
the sheaf on $C$

$$\nb C X ^{\prime}:=Ker(\nb C X \to\hbox{\fascio T}
^1_C)$$

\smallskip\noindent
(see [GK]). The infinitesimal deformations of $C$
into $X$ which preserve the singularity of $C$ are 
parametrized  by $H^0(\nb C X ^{\prime});$ moreover, we have 
$\nb C X ^{\prime}\subset n_{*}\hbox{\fascio N}_f$ (for 
all these facts on $\nb C X ^{\prime}$ see [GK]).
Therefore, $dim(\psf )= 3$ implies

$$3\leqslant h^0(\nb C X ^{\prime})\leqslant h^0(\hbox{
\fascio N}_f),$$

\smallskip\noindent  hence
$h^0(\cO_{\Bbb P^1}(b_1)\oplus\cO_{\Bbb P^1}(b_2)
\geqslant 3.$ From this inequality it follows at once that
the case $b_1<0$ and $b_2<0$ is impossible. 
If $b_1\geqslant 0$ and $b_2\geqslant 0,$ then

$$b_1+b_2+2=h^0(
\cO_{\Bbb P^1}(b_1)\oplus\cO_{\Bbb P^1}(b_2)
)\geqslant 3$$

\smallskip\noindent so $b_1+b_2\geqslant 1$ and
$C^2\geqslant 0.$ Finally, assume $b_1>0$ and
$b_2<0.$ Then $h^0(
\cO_{\Bbb P^1}(b_1)\oplus\cO_{\Bbb P^1}(b_2))
=b_1+1$ and $b_1\geqslant 2.$ 
Assume  $a_1\geqslant b_1,$ so $a_1=b_1+k$ with 
$k\geqslant 0.$ Then, from $(6)$ we get
$b_2+2=a_2+k,$ hence $b_2\geqslant -2.$ From
$b_1\geqslant 2$ it follows $b_1+b_2\geqslant 0,$
and, finally, $C^2\geqslant -1$ by $(4).$

\noindent
So we can assume $a_1< b_1$ which implies $F_{11}=0.$
This forces $F_{12}\neq 0$ and $F_{21}\neq 0, $
hence $a_1-b_2\geqslant 0$ and $a_2-b_1\geqslant 0.$
Moreover, $a_1-b_2>0$ since $b_2<0.$ From $(6)$
it follows that there are only two possibilities:

$$a_1-b_2=1\ \ \ \ \ \ \hbox{and}
\ \ \ \ \ \  a_2-b_1=1$$

\vskip -0.3 cm

$$a_1-b_2=2\ \ \ \ \ \ \hbox{and}
\ \ \ \ \ \  a_2-b_1=0$$

\smallskip\noindent
In the former case $b_2\geqslant -1$ and $b_1+b_2
\geqslant 1,$ \ \ $C^2\geqslant 0.$ In the latter
case we get $C^2\geqslant -1.$

\smallskip\noindent{\tito Step 6.}\ Finally,
we prove that $C^2\geqslant -1$ in the case when
the general $C\in\psf$ is elliptic. First of all, we want 
to show that in this case $\nb C X $ splits. 

\noindent
It is well known that on an elliptic curve $C$ any indecomposable
rank $2$ vector bundle {\fascio E} fits into a short exact
sequence of one of the following forms

$$0\to\hbox{\fascio L} \to\hbox{\fascio E}\to\hbox{\fascio L}\to 0
\ \ \ \ \ \ \ \ \ \ \ \ \
0\to\hbox{\fascio L} \to\hbox{\fascio E}\to\hbox{\fascio L}(P)\to 0$$

\smallskip

\noindent
where {\fascio L} is an invertible sheaf on $C$ and $P\in C.$
The exact sequence $(1)$ becomes

\vskip -0.3 cm

$$0\to\cO_C (C) \to\nb C X\to\cO_C(1)\to 0
\leqno(7)$$

\smallskip

\noindent Assume $\nb C X$ indecomposable. Then $\cO_C (C)
\simeq\hbox{\fascio L}.$ But from $(7)$ it follows at
once that $\cc (\hbox{\fascio L}) \geqslant 2,$ whereas it
was already proved that $\cc (\cO_C (C)) \leqslant 0,$ a
contradiction.

\noindent Therefore, we can
assume \ $\nb C X \simeq 
\hbox{\fascio L}_1\oplus\hbox{\fascio L}_2\ $ 
and set
$d_i:=\cc (\hbox{\fascio L}_i)\ .$ 
Since $\nb C X$ is spanned, both 
$\hbox{\fascio L}_i$
are spanned and $d_i\geqslant 0\ .$ 
If $d_i>0$ for $i=1,2\ ,$ then

\vskip -0.2 cm

$$
C^2+3=\cc (\nb C X )=\cc (\hbox{\fascio L}_1)
+\cc (\hbox{\fascio L}_2)=
h^0(\hbox{\fascio L}_1)+h^0(\hbox{\fascio L}_2)
=h^0(\nb C X )\geqslant 3\ ,$$

\noindent and we can conclude as in the 
previous case.

\smallskip\noindent
Finally, assume \ $\nb C X \simeq 
\cO_C\oplus\hbox{\fascio L}\ $ with 
$\cc (\hbox{\fascio L})>0\ .$ Then, from

$$C^2+3=\cc (\nb C X )=\cc (\hbox{\fascio L})
=h^0(\hbox{\fascio L})$$

\noindent we get

$$C^2+4=h^0(\cO_C)+h^0(\hbox{\fascio L})
=h^0(\nb C X )\geqslant 3$$

\noindent hence $C^2\geqslant -1\ .$
\endproof

\medskip

\noindent{\bf Corollary 3.2}\ {\sl  If $X$ 
satisfies the assumptions of Theorem 1.6,
then $X$ is uniruled.}

\proof
The above lemma jointly with Lemma 2.4 (ii) 
yield $K_X\cdot C<0\ .$ 
Since the curves $C$ cover $X,$ the assertion 
follows from [MM], Thm.1.
\endproof

\medskip

\noindent{\bf Lemma 3.3}\ {\sl 
If\ $K_X+H$ is big and nef, then\ $h^0(\cO_X(K_X+H))
\geqslant 2\ .$ In particular, \ $K_X+H$ is effective
and\ $dim\ \phi_{K_X+H} (X)\geqslant1\ .$}

\proof  By Kodaira's Vanishing Theorem we have 

$$\chi (\cO_X(K_X+H))=h^0(\cO_X(K_X+H))\ .$$

\smallskip

\noindent
Moreover, Serre's duality yields 
$\chi (K_X)=-\chi (\cO_X)\ .$ 
The following inequality, due to Sommese 
([aS], Theorem (1.0)),

$$(K_X+H)^3\leqslant 36 (\chi (K_X)+ 
\chi (\cO_X(K_X+H)))$$

\noindent then becomes

$$h^0(\cO_X(K_X+H))\geqslant {1\over 36} 
(K_X+H)^3 + \chi (\cO_X)\ .$$

\noindent Being $K_X+H$ big and nef, 
$(K_X+H)^3>0\ .$ Therefore, the desired conclusion 
follows since $X$ is uniruled and with 
irregularity zero, which imply
\ $\chi (\cO_X)=1+h^2(\cO_X)\geqslant 1\ .$ 
\endproof

\medskip

We can prove, now, the main result of this
section.

\smallskip 

\noindent{\bf Proposition 3.4}\ {\sl  If $X$ 
satisfies the assumptions of Theorem 1.6,
then $X$ is of log-special type.}

\proof
First of all, we deal with the case $\delta =2.$
From Lemma 2.4,(ii) we get

$$(K_X+H)\cdot C=K_X\cdot C+H\cdot C=-C^2-2$$

\smallskip

\noindent
Then, from Lemma 3.1 it follows at once that
$(K_X+H)$ cannot be nef.

\smallskip

Assume, now, $\delta =3$ and $C^2=0.$ If $K_X+H$
is big and nef, then for some integer $r>>0$ the
rational map 
$\diagram \phi_{r{\scriptscriptstyle {(K_X+H)}}}:X
\rdotted|>\tip & \Bbb P^n \enddiagram$ is 
generically finite onto its image $W.$ Let $D_1,$
$D_2$ be two general (distinct) elements of the
linear system $\vert K_X+H\vert .$ We have $rD_i=
\phi_{r{\scriptscriptstyle {(K_X+H)}}}^{*}L_i,$ 
where $L_i$ is a hyperplane
section of $W,$ for $i=1,2.$ On $W$ we have the 
curve $L_1\cap L_2$ and, since 
$\phi_{r{\scriptscriptstyle {(K_X+H)}}}$
is generically finite over $W,$ then

\vskip -0.2 cm

$$
supp(D_1)\cap supp(D_2)=supp(rD_1)\cap supp(rD_2)=
\phi_{r{\scriptscriptstyle {(K_X+H)}}}^{-1}
(L_1\cap L_2)
$$

\smallskip

\noindent
is also a curve outside the base locus $B$ of 
$\vert K_X+H\vert .$ Therefore, only finitely many
curves $C\in\hbox{\fascio F}$ can be irreducible
components of $(D_1\cap D_2)\backslash B.$ Now, from
$\delta =3$ and from Lemma 2.4 it follows 

\vskip -0.2 cm

$$
(K_X+H)\cdot C=K_X\cdot C+H\cdot C=-C^2
\leqno(8)$$

\smallskip

\noindent
Then, assume $P\in (D_1\cap D_2)\backslash B.$ 
By $(8)$ and the assumption $C^2=0$,  
we have that every curve $C\in\hbox{\fascio F}$
such that $P\in C$ is completely contained both in 
$D_1$ and in $D_2.$ Since we have infinitely many
curves of {\fascio F} containing $P,$ we have
a contradiction.

\smallskip

Finally, assume $C^2=-1,$ always with  $\delta =3.$
By Lemma 3.4 we have a rational map $\diagram 
\phi_{\scriptscriptstyle {K_X+H}}:X\rdotted|>\tip & 
\Bbb P^n \enddiagram\ ,$ for some $n\geqslant 1.$
Let $C$ be a general element of $\psf .$ The divisor
$E$ on $C$ which defines $\phi_{\scriptscriptstyle 
{K_X+H}}\vert_{\scriptscriptstyle C}$ has degree 
$1$ by $(8).$ Therefore, \ 
$h^0(\cO_C(E))=1\ $ and the map
$\phi_{\scriptscriptstyle {K_X+H}}
\vert_{\scriptscriptstyle C}$ is {\sl constant.}
Let $C_{\scriptscriptstyle 0}\in\psf$ be a fixed, 
general curve.  We set

$$\psf^{\thinspace\prime}:=
\{ C\in\psf\ \vert\ C\cap C_{\scriptscriptstyle 0}
\neq\emptyset \}$$

\smallskip

\noindent
$\psf^{\thinspace\prime}$ is a $2$-dimensional 
subfamily of $\psf\ .$

\vskip 0.1 cm

\noindent
We claim that \ $V:=\bigcup_{C\in
\psf^{\thinspace\prime}}C$ {\sl is dense in} 
$X\ .$

\vskip 0.1 cm

\noindent
Assume the contrary. Then $V$ is a surface 
containing a $2$-dimensional family of plane 
cubic curves. By Corrado Segre's theorem ([cS]), 
$V\subset\p\ .$ 
Then, the curves $C\in\psf^{\thinspace\prime}$
belong to the linear system of the plane
sections of $V.$ In fact, otherwise the
degree of $V$ would be at least $4.$ 
Since $X$ is fibered by surfaces like $V,$ 
this would imply the existence on $X$ of a 
family of dimension $\geqslant 3$ of plane 
curves of degree $>3.$ But $X$ is of isolated 
type and this contradicts the bound on the
degree (see remarks after Thm.\thinspace 1.3.)
Therefore, on $V$ we have the $3$-dimensional 
family of its plane sections, and the plane 
 cubic curves on $X$ form a family of
dimension $4,$ which contradicts Lemma 3.1. 
The proof of the claim is complete.

\smallskip

From the claim it follows that the rational map
$\phi_{\scriptscriptstyle {K_X+H}}$ is constant, 
which contradicts $n\geqslant 1.$
\endproof

\bigskip


\tit
{\bf  4.-\ Conics contained in the threefolds 
(E1),$\ldots$,(E5).}
\tit


\medskip

We recall that, if $C\in\hbox{\fascio F}$ is
a conic on a general hyperplane section $S$ 
of $X,$ then $C^2=-1$ by Lemma 3.1, hence $C$ is 
an exceptional curve on $S$.

\medskip

In the case $X$ is a Del Pezzo fibration over 
$\Bbb P^1$ the general hyperplane section $S$ 
of $X$ is an elliptic, regular surface which is
minimal ([cO1]). Therefore, on $S$ we cannot have conics
$C$ with $C^2=-1.$         

\medskip

In both cases (E1),(E2) the generic hyperplane 
section $S$ of $X$ is a non-minimal K3 surface.
Moreover, in both cases the $(-1)$-lines on $S$
are contracted by the adjunction map
$\phi_{\scriptscriptstyle {K_S+H_S}}$ on $S$ 
and  the image $W$ of this map is a 
minimal K3 surface ([AR] and [cO2]). Let us
remark that, by the \lq\lq adjunction 
formula", the restriction to $S$
of the adjunction map
$\phi_{\scriptscriptstyle {K_X+2H}}$ is just
$\phi_{\scriptscriptstyle {K_S+H_S}}.$  Then, 
the image of a conic $C\in\hbox{\fascio F}$ in 
$\phi_{\scriptscriptstyle {K_S+H_S}}$ has degree

$$
C\cdot (K_X+2H)=-C^2-4+4=1$$

\smallskip

\noindent
hence it is a line on $W$ such that
$\phi_{\scriptscriptstyle {K_S+H_S}}(C)^2=-1,$
a contradiction.       
\medskip

Finally, we deal with the cases when $X$ is 
a conic bundle (see [BOSS]). The image $W$ of the map
$\phi_{\scriptscriptstyle {K_X+H}}$ is 
$\Bbb P^2$ in the case $d=9,$ 
while $W\subset\p$ is
a quartic surface in the case $d=12.$ Let
$C\in\hbox{\fascio F}$ be a conic on $X.$
Let $D\in\vert K_X+H\vert ,$ and assume
$D\cap C\neq\emptyset ,$ or, more
precisely, that $P\in D\cap C.$ Then, from
$(K_X+H)\cdot C<0,$ it follows that $C
\subset D.$  Let $R\subset W$ be the
image of $D$ in
$\phi_{\scriptscriptstyle {K_X+H}}.$ 
Moreover, let 
$R^{\prime} \subset W$ be a hyperplane
section different from $R$ and 
containing the point $Q:=
\phi_{\scriptscriptstyle {K_X+H}}(P).$
Finally, set $D^{\prime}\in\vert K_X+H\vert$
the preimage of $R^{\prime}.$ 
Then, clearly, $C\subseteq D\cap D^{\prime},$
hence $C$ is a conic of the fibration.
Therefore, $dim (\psf )=2.$

\medskip

We can conclude that in each case (E1), 
$\ldots$ , (E5) we cannot have on $X$ a
family of conics of dimension $3.$ 

\bigskip


\tit
{\bf  5.-\ Plane cubics contained in the threefolds 
(E1),$\ldots $,(E5). }
\tit


\medskip

In this section we will prove that a $3$-fold of one of the types
(E1),$\dots $,(E5) cannot contain a family of dimension $3$
of plane cubics.
\smallskip
\item{(E1)} To rule out this $3$-fold, we need to 
 study something more
its geometry. 

 It is well known (see e.g.  [C] and [AR]) that such a $X$
 is ruled by lines over a surface 
$Y\simeq \Bbb G(1,5)\cap\Bbb P^8$, 
and precisely
$X\simeq \Bbb P(U^*\mid _Y)$ where $U$ 
is the universal bundle on 
the Grassmannian. A general hyperplane
 section $S$ is a non-minimal K3 surface
 with $\pi= 8$ and $K^2_S=-5$. Therefore
 $S$ contains five
 $(-1)$-lines $E_1$,...,$E_5$ which 
are blown-down by the adjunction
 map $\phi :=\phi_{\scriptscriptstyle {K_S+H_S}}:
S\longrightarrow \Bbb P^8$ (note that
 $K_S =E_1+...+E_5$); 
the image of $\phi$, which is a
 minimal K3 surface $S'$ of degree
 $14$
 with $\pi =8$, coincides with 
the previous $Y$. 

 Let us assume now that a K3 
scroll $X$ contains 
a 3-dimensional family of plane 
cubics: then $S$ contains
 at least one cubic $C$
 of the family. Let $C':=\phi 
(C)\subset S'.$ In the case $C^2=0$, 
if $H'$ is a general 
hyperplane
 of $\Bbb P^8$, we have $$C'\cdot H' = 
(\phi ^*C')\cdot(\phi ^*H') = 
(\phi ^*C')\cdot(H+K_S)=C\cdot H+\sum 
(\phi ^*C')\cdot E_i = C\cdot H.$$
 Hence also $C'$ is a plane cubic,
 birational to $C$.
 Let $\Gamma '$ be a hyperplane section of 
$S'$:  $\Gamma '$ is a canonical curve
which possesses at least one 
trisecant line. 
It is easy to see that if the
 canonical curve $\Gamma '$ has one
 trisecant line, 
then it has a 1-dimensional
 family of trisecants, that 
cut on $\Gamma '$ a $g^1_3$; 
so $\Gamma '$ is trigonal and its 
trisecants generate a rational normal scroll.
In this case the  homogeneous 
ideal of $S'$
 is generated by quadrics and 
cubics and the quadrics 
containing $S'$ intersect along 
a rational normal scroll $V$ 
(possibly a cone), such that 
the fibers of the
 restriction to $S'$ of the
 natural map $V\rightarrow
\Bbb P^1$ are precisely the curves of
$\mid C'\mid$. The diagram 
$$
\objectmargin {0.4pc}
\diagram
S\rto^\phi &S'\drto\rmapsin &V \dto\cr 
&&\Bbb P^1
\enddiagram
$$

\noindent shows that also $S$ is
ruled by plane cubics. So, by [kR],
 $S$ is  contained in a
 quadric: but this is  impossible by [AR] (2.12).

In the case $C^2=-1$, by \lq\lq adjunction formula" 
$C\cdot K=1,$ namely $C$ intersects exactly 
one exceptional line. From this it follows that 
$\mid C'\mid$ is a pencil of elliptic quartics.
This pencil cuts on the general hyperplane section
$\Gamma '$ of $S'$ a $g^1_4$, so $\Gamma '$ is a 
canonical tetragonal curve. By [sM], it follows that 
 $\Gamma '$ is not a linear section of 
$\Bbb G(1,5)$. This contradiction 
shows that also this case is impossible. 
\medskip

\item{(E2)}  $X$ has degree $d=7$ and sectional genus
$\pi =5.$

The adjunction map $\phi_{K+2H}:X\to 
\Bbb P^6$
is the contraction of the exceptional 
plane contained in $X$, and the image $X'$
is a Fano manifold, complete intersection of type $(2,2,2)$.
 As in the previous case, 
if $X$ contains a $3$-dimensional 
family of plane cubics and 
if $C^2=0$, then also $X'$ contains  a similar family. 
Hence a general hyperplane section
$S'$ of $X'$ contains at least one plane 
cubic and a general hyperplane section $\Gamma '$ of $S'$
 possesses at least
one trisecant line. But $\Gamma '$ is a canonical curve,
so it is trigonal. From the classification 
of trigonal Fano varieties ([I]), it follows that 
this case is impossible.

If $C^2=-1$, as in the previous case we may conclude that 
$S'$ contains a pencil of elliptic quartics curves.
The union of the $\p$'s generated by such quartics is a
quadric of rank $4$ in $\ppp$ containing $S'$. Since $S'$ is 
a general hyperplane section of $X'$, also $X'$ 
is contained in a quadric $Q$ of rank $4$. Let $X'=
Q\cap Q_1\cap Q_2$; then an elementary computation
 shows that $X'$ is singular at the four points
of intersection of the plane $Sing Q$ with $ Q_1\cap Q_2$:
also this case is excluded.

\medskip

\item{(E3)} $X$ is a Del Pezzo fibration over $\Bbb P^1.$

It is shown in [BOSS], that the 
image of the map  $\phi_{K+H}$
is $\Bbb P^1$, hence the case $C^2=-1$ is excluded as in 
the proof of Prop.3.5. 

Moreover the fibers 
 are Del Pezzo surfaces of $\pp$, i. e. complete intersections
of type $(2,2)$.  
If  $C^2=0$ then, from
 $(K+H)\cdot C\leqslant 0$, it follows that the plane 
cubics should be contained in such Del Pezzo surfaces, so each 
surface of $\mid K+H\mid$ should contain a family of
dimension $2$ of plane cubics: but this contradicts 
the theorem of 
C.Segre ([cS]).

\medskip

\item{(E4),}(E5) If $C^2=-1$, we may argue as in the proof of
Prop.3.5. If $C^2=0$, then the discussion goes as in
 the case of the family of conics.

\vfill\eject

\centerline{\bf References}

\medskip

\noindent\item {[AR]}   Aure A. B. -
 Ranestad K.: 
 The smooth surfaces of degree $9$ 
in $\pp$.
Complex Projective Geometry, Proceedings
 Bergen-Trieste, London Math. Soc. LNS {\bf 179}, 
32-46  (1992)

\vskip 0.1 cm 

\noindent\item {[BOSS]} Braun R. - Ottaviani G. - Schneider M.
- Schreyer F.O.: 
 Classification of log-special $3$-folds in $\ppp$.
Preprint Bayreuth (1992)

\vskip 0.1 cm 

\noindent\item {[BH]}  Bruns W. - Herzog J.: 
Cohen-Macaulay rings.
Cambridge Univ. Press, 1993

\vskip 0.1 cm

\noindent\item {[C]}  Chang M.C.:
Classification of Buchsbaum subvarieties of 
codimension $2$ in projective space.
J. reine angew. Math. {\bf 401},  
101-112 (1989)

\vskip 0.1 cm

\noindent\item {[DP]}  Decker W. - Popescu S.:
On surfaces in $\Bbb P^4$ and $3$-folds in $\Bbb P^5$. 
Preprint Saarbr\"ucken (1993)

\vskip 0.1 cm

\noindent\item {[aF]}  Franchetta A.:
Sulla curva doppia
della proiezione della superficie 
generale dell'$S_{\scriptstyle
4},$ da un punto generico su un 
$S_{\scriptstyle 3}$.
Rend. Acc. d'Italia ser. VII-2,  282-288 (1940)

\vskip 0.1 cm 

\noindent\item {[tF]}   Fujita T.:\ 
Classification Theories of Polarized Varieties.
London Math. Soc. Lecture Note Ser.
155, Cambridge University Press, 1990

\vskip 0.1 cm

\noindent\item {[GK]} Greuel G.M. -  Karras U.:
Families of varieties with prescribed singularities.
Comp. Math.  {\bf 69}, 83-110 (1989)

\vskip 0.1 cm

\noindent\item {[I]}  Iskovskih V.A.:
Fano $3$-folds II.
Math. U.S.S.R. Izv. {\bf 12}, 469-506 (1978)

\vskip 0.1 cm

\noindent\item {[LT]}  Lanteri A. - Turrini C.:
Some formulas concerning non-singular
 algebraic varieties embedded in some ambient variety.
Atti Accad. Naz. Lincei (8) {\bf 69},
 463-474 (1980)

\vskip 0.1 cm

\noindent\item {[eM]}  Mezzetti E.:
Projective \  varieties \ with many
degenerate \  subvarieties.
BUMI (7) {\bf 8-B},  807-832 (1994)

\vskip 0.1 cm

\noindent\item {[MP]}  Mezzetti E. -  Portelli D.: 
 A tour through some classical theorem on algebraic
surfaces. Preprint, Trieste 1995

\vskip 0.1 cm

\noindent\item {[MR]}  Mezzetti E. -  Raspanti I.:
A Laudal-type theorem for surfaces in $\Bbb P^4$.
Rend. Sem. Mat. Univ. Politec. Torino 
{\bf 48}, 529-537 (1990)

\vskip 0.1 cm

\noindent\item {[MM]}   Miyaoka Y. -  Mori S.:
A numerical criterion for uniruledness.
Annals of Math. {\bf 124}, 65-69 (1986)

\vskip 0.1 cm

\noindent\item {[sM]}   Mukai S.:
Curves and Symmetric Spaces.
Proc. Japan Acad. {\bf 68} Ser.A, 
 7-10 (1992)

\vskip 0.1 cm

\noindent\item {[cO1]}  Okonek C.:
 Fl\"achen vom Grad $8$ in $\pp$.
Math Z. {\bf 191},  207-223 (1986)

\vskip 0.1 cm

\noindent\item {[cO2]}  Okonek C.:
\"Uber $2$-codimensionale Untermannigfaltigkeiten
vom Grad $7$ in $\pp$ und $\ppp$.
Math Z. {\bf 187},  209-219 (1984)

\vskip 0.1 cm

\noindent\item {[gO]}  Ottaviani G.:
 On $3$-folds in $\ppp$ which are scrolls.
Ann. Scuola Norm. Sup. Pisa Cl.
Sci. Ser. (4) {\bf 19},  451-471 (1992)

\vskip 0.1 cm

\noindent\item {[kR]} Ranestad K.:
On Smooth Plane Curve Fibrations in $\pp$.
In Geometry of Complex Projective Varieties,
Cetraro (Italy) 1990, Mediterranean Press 
 243-255 (1993)

\vskip 0.1 cm

\noindent\item {[eR]} Rogora E.:
Varieties with many lines.
Manuscr. Math. {\bf 82},  207-226 (1994)

\vskip 0.1 cm

\noindent\item {[bS]}  Segre B.:
 Sulle $V_n$ contenenti pi\'u
di $\infty ^{n-k} S_k$.  
I, Lincei - Rend. Sc. fis.
mat. e nat. {\bf 5}, 193-197 (1948);  
II, Lincei - Rend. Sc. fis.
mat. e nat. {\bf 5}, 275-280 (1948)

\vskip 0.1 cm
\noindent\item {[cS]}    Segre C.: 
Le superficie degli
iperspazi con una doppia infinit\` a 
di curve piane o spaziali.  
Atti R. Acc. 
Scienze di Torino {\bf 56},
 75-89 (1921)

\vskip 0.1 cm

\noindent\item {[fS]}  Severi F.:
 Intorno ai punti doppi
improprii di una superficie 
generale dello spazio a quattro 
dimensioni, e ai suoi punti tripli 
apparenti. Rend. Circolo 
Matematico di Palermo {\bf 15},
 33-51 (1901)

\vskip 0.1 cm

\noindent\item {[aS]}  Sommese A. J.: 
On the nonemptiness of the adjoint
linear system of a hyperplane section of a threefold.
J. reine angew. Math. {\bf 402},  211-220 (1989)

\vskip 0.1 cm

\noindent\item {[VdV]} Van de Ven  A.: 
On the $2$-connectedness of a 
very ample divisor on a surface.
Duke Math. J.
{\bf 46},  403-407 (1979)

\bye